\documentclass[a4paper,12pt]{amsart}
\usepackage{amssymb}
\usepackage{pifont} 
\usepackage{bbding}
\usepackage{wasysym}
\usepackage{ifthen}
\usepackage{cite}
 \usepackage[dvips]{graphicx}
\nonstopmode \numberwithin{equation}{section}
\usepackage{amssymb}
\usepackage{ifthen}
\usepackage{graphicx}
\usepackage{amsmath}
\usepackage[T1]{fontenc} 
\usepackage[utf8]{inputenc}
\usepackage[usenames,dvipsnames]{color}
\usepackage{color}
\usepackage[english]{babel}
\usepackage{fancyhdr}
\usepackage{fancybox}
\usepackage{tikz}

\nonstopmode \numberwithin{equation}{section}
\theoremstyle{plain}

\newtheorem{conj}{Conjecture}

\theoremstyle{definition}

\newtheorem{example}{Example}[section]
\newtheorem{thm}{Theorem}[section]
\newtheorem{prob}{Problem}[section]
\newtheorem{cor}{Corollary}[section]
\newtheorem{ques}{Question}[section]
\newtheorem{prop}{Proposition}[section]
\newtheorem{rem}{Remark}[section]
\newtheorem{lem}{Lemma}[section]

\newcounter{minutes}
\divide\time by 60
\newcounter{hours}
\multiply\time by 60
\addtocounter{minutes}{-\time}

\newcounter {own}
\def\theown {\thesection       .\arabic{own}}

\newenvironment{pf}[1][]{%
 \vskip 3mm
 \noindent
 \ifthenelse{\equal{#1}{}}%
  {{\slshape Proof. }}%
  {{\slshape #1.} }%
 }%
{\qed\bigskip}

\newcounter{alphabet}

\def\be{\begin{equation}}
\def\ee{\end{equation}}

\newcommand{\bee}{\begin{enumerate}}
\newcommand{\eee}{\end{enumerate}}

\newcommand{\blem}{\begin{lem}}
\newcommand{\elem}{\end{lem}}
\newcommand{\bthm}{\begin{thm}}
\newcommand{\ethm}{\end{thm}}
\newcommand{\bcor}{\begin{cor}}
\newcommand{\ecor}{\end{cor}}
\newcommand{\beg}{\begin{examp}}
\newcommand{\eeg}{\end{examp}}
\newcommand{\begs}{\begin{examples}}
\newcommand{\eegs}{\end{examples}}

\newcommand{\bdefn}{\begin{defn}}
\newcommand{\edefn}{\end{defn}}

\newcommand{\bprob}{\begin{prob}}
\newcommand{\eprob}{\end{prob}}
\newcommand{\bei}{\begin{itemize}}
\newcommand{\eei}{\end{itemize}}

\newcommand{\bcon}{\begin{conj}}
\newcommand{\econ}{\end{conj}}
\newcommand{\bcons}{\begin{conjs}}
\newcommand{\econs}{\end{conjs}}
\newcommand{\bprop}{\begin{prop}}
\newcommand{\eprop}{\end{prop}}
\newcommand{\br}{\begin{rem}}
\newcommand{\er}{\end{rem}}
\newcommand{\brs}{\begin{rems}}
\newcommand{\ers}{\end{rems}}
\newcommand{\bo}{\begin{obser}}
\newcommand{\eo}{\end{obser}}
\newcommand{\bos}{\begin{obsers}}
\newcommand{\eos}{\end{obsers}}
\newcommand{\bpf}{\begin{pf}}
\newcommand{\epf}{\end{pf}}
\newcommand{\ba}{\begin{array}}
\newcommand{\ea}{\end{array}}
\newcommand{\beq}{\begin{eqnarray}}
\newcommand{\beqq}{\begin{eqnarray*}}
\newcommand{\eeq}{\end{eqnarray}}
\newcommand{\eeqq}{\end{eqnarray*}}

\begin{document}

\title{Hankel and Toeplitz determinants of logarithmic coefficients of Inverse functions for certain classes of univalent functions}

\author{Sanju Mandal}
\address{Sanju Mandal, Department of Mathematics, Jadavpur University, Kolkata-700032, West Bengal,India.}
\email{sanjum.math.rs@jadavpuruniversity.in, sanju.math.rs@gmail.com}

\author{Partha Pratim Roy}
\address{Partha Pratim Roy, Department of Mathematics, Jadavpur University, Kolkata-700032, West Bengal, India.}
\email{pproy.math.rs@jadavpuruniversity.in}

\author{Molla Basir Ahamed}
\address{Molla Basir Ahamed, Department of Mathematics, Jadavpur University, Kolkata-700032, West Bengal,India.}
\email{mbahamed.math@jadavpuruniversity.in}

\subjclass[{AMS} Subject Classification:]{Primary 30A10, 30H05, 30C35, Secondary 30C45}

\def\thefootnote{}
\footnotetext{ {\tiny File:~\jobname.tex,
printed: \number\year-\number\month-\number\day,
          \thehours.\ifnum\theminutes<10{0}\fi\theminutes }
} \makeatletter\def\thefootnote{\@arabic\c@footnote}\makeatother
\begin{abstract} 
The Hankel and Toeplitz determinants $H_{2,1}(F_{f^{-1}}/2)$ and $T_{2,1}(F_{f^{-1}}/2)$ are defined as:
\begin{align*}
	H_{2,1}(F_{f^{-1}}/2):= \begin{vmatrix}
		\Gamma_1 & \Gamma_2 \\
		\Gamma_2 & \Gamma_3
	\end{vmatrix} 
	\;\;\mbox{and} \;\; 
	T_{2,1}(F_{f^{-1}}/2):= \begin{vmatrix}
		\Gamma_1 & \Gamma_2 \\
		\Gamma_2 & \Gamma_1
	\end{vmatrix}
\end{align*}
where $\Gamma_1, \Gamma_2,$ and $\Gamma_3$ are the first, second and third logarithmic coefficients of inverse functions belonging to the class $\mathcal{S}$ of normalized univalent functions. In this article, we establish sharp inequalities $|H_{2,1}(F_{f^{-1}}/2)|\leq 1/4$, $|H_{2,1}(F_{f^{-1}}/2)| \leq 1/36$, $|T_{2,1}(F_{f^{-1}}/2)|\leq 5/16$  and $|T_{2,1}(F_{f^{-1}}/2)|\leq 145/2304$ for the logarithmic coefficients of inverse functions for the classes starlike functions and convex functions with respect to symmetric points. In addition, our findings are substantiated further through the incorporation of illustrative examples, which support the strict inequality and lend credence to our conclusions.
\end{abstract}
\maketitle
\pagestyle{myheadings}
\markboth{S. Mandal, P. P. Roy,  and  M. B. Ahamed}{Hankel and Toeplitz determinants of logarithmic coefficients of Inverse functions}
\tableofcontents
\section{Introduction}
Let $\mathcal{H}$ be the class of functions $ f $ which are holomorphic in the open unit disk $\mathbb{D}=\{z\in\mathbb{C}: |z|<1\}$ of the form 
\begin{align}\label{eq-1.1}
	f(z)=\sum_{n=1}^{\infty}a_nz^n,\; \mbox{for}\; z\in\mathbb{D}.
\end{align}
Then $\mathcal{H}$ is a locally convex topological vector space endowed with the topology of uniform convergence over compact subsets of $\mathbb{D}$. Let $\mathcal{A}$ denote the class of function $f\in\mathcal{H}$ such that $f(0)=0$ and $f^{\prime}(0)=1$. That is, the function $f$ takes the following  form
\begin{align}\label{eq-1.2}
	f(z)=z+ \sum_{n=2}^{\infty}a_nz^n,\; \mbox{for}\; z\in\mathbb{D}.
\end{align} 
Let $\mathcal{S}$ denote the subclass of all functions in $\mathcal{A}$ which are univalent. For a general theory of univalent functions, we refer the classical books \cite{Duren-1983-NY,Goodman-1983}.

Let
\begin{align}\label{eq-1.3}
	F_{f}(z):=\log\dfrac{f(z)}{z}=2\sum_{n=1}^{\infty}\gamma_{n}(f)z^n, \;\; z\in\mathbb{D},\;\;\log 1:=0,
\end{align}
be a logarithmic function associated with $f\in\mathcal{S}$.  The numbers $\gamma_{n}:=\gamma_{n}(f)$ ($ n\in\mathbb{N} $) are called the logarithmic coefficients of $f$. Although the logarithmic coefficients $\gamma_{n}$ play a critical role in the theory of univalent functions, it appears that only a limited number of sharp bounds have been established for them. As is well known, the logarithmic coefficients play a crucial role in Milin’s conjecture (\cite{Milin-1977-ET}, see also \cite[p.155]{Duren-1983-NY}).
Milin conjectured that for $f\in\mathcal{S}$ and $n\geq2$,
\begin{align*}
	\sum_{m=1}^{n}\sum_{k=1}^{m}\left(k|\gamma_{k}|^2 -\frac{1}{k}\right)\leq 0
\end{align*}
where the equality holds if, and only if, $f$ is a rotation of the Koebe function. De Branges \cite{Branges-AM-1985} proved Milin conjecture which confirmed the famous Bieberbach conjecture. On the other hand, one of reasons for more attention has been given to the logarithmic coefficients is that the sharp bound for the class $\mathcal{S}$ is known only for $\gamma_{1}$ and $\gamma_{2}$, namely
\begin{align*}
	|\gamma_{1}|\leq 1, \;\; |\gamma_{2}|\leq \dfrac{1}{2}+ \dfrac{1}{e} =0.635\ldots
\end{align*}
It is known that for the Koebe function $f(z)=z/(1-z)^2$, the logarithmic coefficients are $\gamma_{n}=1/n$, for each positive integer $n$. Since the Koebe function appears as an extremal function in many problems of geometric theory of analytic functions, one could expect that $\gamma_{n}=1/n$ holds for functions in $\mathcal{S}$. But this is not true in general, even in order of magnitude. Finding the sharp bounds of $\gamma_{n}$, with $n\geq 3$, for the class $\mathcal{S}$ remains an unsolved problem. Estimating the modulus of logarithmic coefficients for $f\in\mathcal{S}$ and various sub-classes has been considered recently by several authors. For in-depth information on the topic, we suggest referring to the articles \cite{Ali-Allu-PAMS-2018,Ali-Allu-Thomas-CRMCS-2018,Cho-Kowalczyk-kwon-Lecko-Sim-RACSAM-2020,Thomas-PAMS-2016,Ponnusamy-Sugawa-BDSM-2021,Roth-PAMS-2007,Girela-AASF-2000}, and exploring the references provided therein.\vspace{2mm}

A function $f$ is said to be starlike in $\mathbb{C}$ if it maps the open unit disk $\mathbb{D}:=\{z\in\mathbb{C} : |z|<1\}$ conformally onto a region that is starlike with respect to the origin. The class of functions that are starlike with respect to symmetric points, which is denoted as $\mathcal{S}^*_S$, was introduced by Sakaguchi in $1959$ \cite{Sakaguchi-JMSJ-1959}. A function $f$ that belongs to $\mathcal{S}^*_S$ is characterized by the following conditions:
\begin{align*}
	\mbox{Re}\left(\dfrac{zf^{\prime}(z)}{f(z)-f(-z)}\right)> 0, \;\;z\in\mathbb{D}.
\end{align*}

A function $f$ is said to be convex in $\mathbb{C}$ if it maps the open unit disk $|z|<1$ conformally onto a region that is convex. We consider another class, which is denoted by $\mathcal{K}_S$, that is, a function $f\in\mathcal{A}$ is said to be convex with respect to symmetric points if, and only if,
\begin{align*}
	\mbox{Re}\left(\dfrac{(zf^{\prime}(z))^{\prime}}{(f(z)-f(-z))^{\prime}}\right)> 0, \;\;z\in\mathbb{D}.
\end{align*}
The functions that are members of $\mathcal{S}^*_S$ are identified as close-to-convex, and consequently, they are univalent. The class of functions that are starlike with respect to symmetric points also includes the classes of convex functions and odd starlike functions with respect to the origin.\vspace{1.2mm} 

The study of coefficients problems in geometric function theory has a rich history. The challenging and interesting task lies in determining sharp bounds for various inequalities related to functions within different classes of univalent functions. Over the past several years, a significant amount of research in geometric function theory has been devoted to establishing sharp bounds of Hankel or Toeplitz determinants for different classes of functions. However, upon in-depth analysis of these findings, it becomes evident that there is relatively less study focused on exploring the aforementioned determinants not only for logarithmic coefficients but also for the inverse functions of these classes. This study is dedicated to providing the sharp bound for the second Hankel and Toeplitz determinants, whose entries are the logarithmic coefficients of inverse functions.\vspace{1.2mm}

Let $F$ be the inverse function of $f\in\mathcal{S}$ defined in a neighborhood of the origin with the Taylor series expansion
\begin{align}\label{eq-1.4}
	F(w):=f^{-1}(w)= w+\sum_{n=2}^{\infty} A_n w^n,
\end{align}
where we may choose $|w|<1/4$, as we know that the famous Koebe’s $1/4$-theorem ensures that, for each univalent function $f$ defined in $\mathbb{D}$, it inverse $f^{-1}$ exists at least on a disc of radius $1/4$. Using a variational method, L$\ddot{\mbox{o}}$wner \cite{Lowner-IMA-1923} obtained the sharp estimate $ |A_n|\leq K_n \;\mbox{for each}\; n, $ where $K_n=(2n)!/(n!(n+1)!)$ and $K(w)= w +K_2 w^2 +K_3 w^3 +\cdots$ is the inverse of the K\"oebe function. There has been a good deal of interest in determining the behavior of the inverse coefficients of $f$ given in \eqref{eq-1.2} when the corresponding function $f$ is restricted to some proper geometric subclasses of $\mathcal{S}$.\vspace{2mm}

Let $f(z)=z+ \sum_{n=2}^{\infty}a_nz^n$ be a function in class $\mathcal{S}$. Science $f(f^{-1}(w))=w$ and using \eqref{eq-1.4} we obtain
\begin{align}\label{eq-1.5}
	\begin{cases}
		A_2= -a_2, \\ A_3=-a_3 +2a^2_{2}, \\ A_4=- a_4 +5a_2 a_3 -5a^3_{2}, \\ A_5=- a_5+6a_4 a_2- 21a_3 a^2_{2} +3a^2_{3} +14a^4_{2}.
	\end{cases}
\end{align}
The notation of the logarithmic coefficient of inverse of $f$ was introduced by Ponnusamy \textit{et al.} (see \cite{Ponnusamy-Sharma-Wirths-RM-2018}). As with $f$, the logarithmic inverse coefficients $\Gamma_n$, $n\in\mathbb{N}$, of $F$ are defined by the equation
\begin{align}\label{eq-1.6}
	\log\left(\frac{F(w)}{w}\right)=2\sum_{n=1}^{\infty} \Gamma_n(F) w^n \; \mbox{for} \;|w|<1/4.
\end{align}
The author's obtained the sharp bound for the logarithmic inverse coefficients for the class $\mathcal{S}$. In fact, Ponnusamy \textit{et al.} \cite{Ponnusamy-Sharma-Wirths-RM-2018} have established for $f\in\mathcal{S}$ that
\begin{align*}
	|\Gamma_n(F)|\leq\frac{1}{2n}\binom{2n}{n}
\end{align*}
and showed that the equality holds only for the K\"oebe function or  its rotations. Differentiating \eqref{eq-1.6} and using \eqref{eq-1.5}, we obtain
\begin{align}\label{eq-1.7}
	\begin{cases}
		\Gamma_1=-\frac{1}{2}a_2, \\ \Gamma_2=-\frac{1}{2}\left(a_3 -\frac{3}{2}a^2_{2}\right), \\ \Gamma_3=-\frac{1}{2}\left(a_4 -4a_2 a_3 +\frac{10}{3}a^3_{2}\right), \\ \Gamma_4 =-\frac{1}{2} \left(a_5 -5a_4 a_2 +15 a_3 a^2_{2} -\frac{5}{2}a^2_{3} -\frac{35}{4}a^4_{2}\right).
	\end{cases}
\end{align}
 In $2016$, Ye and Lim \cite{Ye-Lim-FCM-2016} proved that any $n\times n$ matrix over $\mathbb{C}$ generically can be written as the product of some Toeplitz matrices or Hankel matrices. Hankel matrices and determinants are essential components in various mathematical disciplines, offering a wide range of applications \cite{Ye-Lim-FCM-2016}. Likewise, in both pure and applied mathematics, Toeplitz matrices and Toeplitz determinants play a role \cite{Toeplitz-1907}. They occur in analysis, integral equations, image processing, signal processing, quantum mechanics and among other areas. For more applications, we refer to the survey article \cite{Ye-Lim-FCM-2016} and references therein. Toeplitz matrices contain same entries along their diagonals. This study is dedicated to providing the sharp bound for the second Hankel and Toeplitz determinants, whose entries are the logarithmic coefficients.\vspace{1.2mm}

 The Hankel and Toeplitz determinants $H_{q,n}(f)$ and $T_{q,n}(f)$ of Taylor's coefficients of functions $f\in\mathcal{A}$ represented by \eqref{eq-1.1} is defined for $q,n\in\mathbb{N}$ as follows:
\begin{align*}
	H_{q,n}(f):=\begin{vmatrix}
		a_{n} & a_{n+1} &\cdots& a_{n+q-1}\\ a_{n+1} & a_{n+2} &\cdots& a_{n+q} \\ \vdots & \vdots & \vdots & \vdots \\ a_{n+q-1} & a_{n+q} &\cdots& a_{n+2(q-1)}
	\end{vmatrix},\;
	T_{q,n}(f):=\begin{vmatrix}
		a_{n} & a_{n+1} &\cdots& a_{n+q-1}\\ a_{n+1} & a_{n} &\cdots& a_{n+q-2} \\ \vdots & \vdots & \vdots & \vdots \\ a_{n+q-1} & a_{n+q-2} &\cdots& a_{n}
	\end{vmatrix}.
\end{align*}
Recently, Kowalczyk and Lecko \cite{Kowalczyk-Lecko-BAMS-2022} proposed a Hankel determinant whose elements are the logarithmic coefficients of $f\in\mathcal{S}$, realizing the extensive use of these coefficients. Inspired by these ideas, we introduce the investigation of the Hankel determinant $H_{q,n}(F_{f^{-1}}/2)$ and Toeplitz determinant $T_{q,n}(F_{f^{-1}}/2)$, wherein the elements are logarithmic coefficients of inverse functions of $f^{-1}\in\mathcal{S}$. The determinant $H_{q,n}(F_{f^{-1}}/2)$ is expressed as follows:
\begin{align*}
	H_{q,n}(F_{f^{-1}}/2)=\begin{vmatrix}
		\Gamma_{n} & \Gamma_{n+1} &\cdots& \Gamma_{n+q-1}\\ \Gamma_{n+1} & \Gamma_{n+2} &\cdots& \Gamma_{n+q} \\ \vdots & \vdots & \vdots & \vdots \\ \Gamma_{n+q-1} & \Gamma_{n+q} &\cdots& \Gamma_{n+2(q-1)}
	\end{vmatrix},
\end{align*}
and, the determinant $T_{q,n}(F_{f^{-1}}/2)$ is expressed as follows:
\begin{align*}
	T_{q,n}(F_{f^{-1}}/2)=\begin{vmatrix}
		\Gamma_{n} & \Gamma_{n+1} &\cdots& \Gamma_{n+q-1}\\ \Gamma_{n+1} & \Gamma_{n} &\cdots& \Gamma_{n+q-2} \\ \vdots & \vdots & \vdots & \vdots \\ \Gamma_{n+q-1} & \Gamma_{n+q-2} &\cdots& \Gamma_{n}
	\end{vmatrix}.
\end{align*}
The study of Hankel and Toeplitz determinants for starlike, convex, or many other functions has been done extensively (see \cite{Ponnusamy-Sugawa-BDSM-2021,Kowalczyk-Lecko-RACSAM-2023,Raza-Riza-Thomas-BAMS-2023,Sim-Lecko-Thomas-AMPA-2021,Kowalczyk-Lecko-BAMS-2022, Allu-Lecko-Thomas-MJM-2022,Ali-Thomas-Allu-BAMS-2018,Kowalczyk-Lecko-Sim-BAMS-2018}), their sharp bounds have been established. Recently, Hankel determinants with logarithmic coefficients have been studied for certain subclasses of starlike, convex, univalent, strongly starlike, and strongly convex functions (see \cite{Allu-Arora-Shaji-MJM-2023,Kowalczyk-Lecko-BAMS-2022,Kowalczyk-Lecko-BMMS-2022,Libera-Zlotkiewicz-PAMS-1982} and references therein). Despite the significance of these problems, the sharp bound of Toeplitz determinants with logarithmic coefficients of inverse functions have not been investigated. Some progress has been made in this direction, especially concerning Toeplitz determinants for univalent functions with specific symmetries or in particular special domains. In 2021, Zaprawa \cite{Zaprawa-BSMM-2021} obtained the sharp bounds of the initial logarithmic coefficients $\gamma_{n}$ for functions in the classes $\mathcal{S}^*_S$  and $\mathcal{K}_S$. However, no work has been carried out for sharp bounds of Hankel and Toeplitz determinants of logarithmic coefficients of inverse functions lying in starlike functions and convex functions with respect to symmetric points.\vspace{2mm} 

However, the primary aim of this paper is to find the sharp bound of
\begin{align}\label{eq-1.8}
	H_{2,1}(F_{f^{-1}}/2)=\Gamma_{1}\Gamma_{3} -\Gamma^2_{2}=\frac{1}{48} \left(13a^4_2 -12a^2_2 a_3 - 12 a^2_3 + 12 a_2 a_4\right),
\end{align} 
and 
\begin{align}\label{eq-1.9}
	T_{2,1}(F_{f^{-1}}/2)=\Gamma^2_{1}-\Gamma^2_{2}= \frac{1}{16}(-9a^4_2 + 4a^2_2 -4a^2_3 + 12a^2_2 a_3)
\end{align}
when $f^{-1}$ is a member of the classes $\mathcal{S}^*_S$ and $\mathcal{K}_S$, which respectively refer to starlike functions and convex functions with respect to symmetric points. In addition, we give examples of functions to illustrate these results.

\section{Main results for classes $\mathcal{K}_S$ and $ \mathcal{S}^{*}_S $, remarks and examples}
Suppose $\mathcal{B}_0$ be the class of Schwarz function \textit{i.e} analytic function $w:\mathbb{D}\rightarrow\mathbb{D}$ such that $w(0)=0$. A function $w\in\mathcal{B}_0$ can be written as a power series
\begin{align*}
	w(z)=\sum_{n=1}^{\infty}c_n z^n.
\end{align*}
For two functions $f$ and $g$ which are analytic in $\mathbb{D}$, we say that the function $f$ is subordinate to $g$ in $\mathbb{D}$ and written as
$f(z)\prec g(z)$ there exists a Schwarz function $w\in\mathcal{B}_0$ such that $f(z)=g(w(z))$, $z\in\mathbb{D}$. To prove our results, the following lemma for Schwarz functions will play a key role.
\begin{lem}\cite{Efraimidis-JMAA-2016}\label{lem-2.1}
Let $ w(z)=c_1 z +c_2 z^2 + c_3 z^3+\ldots $ be a Schwarz function. Then 	
\begin{align*}
	|c_1|\leq 1, \;\;|c_2|\leq 1 -|c_1|^2, \;\;\mbox{and}\;\; |c_3|\leq 1 -|c_1|^2 -\dfrac{|c_2|^2}{1+|c_1|}.
\end{align*}
\end{lem}
\subsection{Sharp bound of $|H_{2,1}(F_{f^{-1}}/2)|$ for $f\in \mathcal{S}^*_S$:}
We obtain the following sharp bound for $H_{2,1}(F_{f^{-1}}/2)$ for functions in the class $ \mathcal{S}^{*}_S $.
\begin{thm}\label{th-2.1}
Let $ f\in \mathcal{S}^{*}_S $. Then
\begin{align*}
	|H_{2,1}(F_{f^{-1}}/2)|\leq\dfrac{1}{4}.
\end{align*}
The inequality is sharp with the extremal function	
\begin{align*}
	h_1(z)=\frac{z}{1-z^2}= z +z^3 +z^5 +\cdots, \;\;\;\;z\in\mathbb{D}.
\end{align*}
\end{thm}
An example is presented below to demonstrate that the strict inequality in Theorem \ref{th-2.1} remains valid.
\begin{example}
	Consider the function
	\begin{align*}
		h_2(z)=\frac{z}{1 -z}=z +z^2 +z^3 +\cdots, \;\;\;\;z\in\mathbb{D}.
	\end{align*}
	It is easy to see that 
	\begin{align*}
		\mbox{Re}\left(\dfrac{zh^{\prime}_{2}(z)}{h_{2}(z)-h_{2}(-z)}\right) =\dfrac{1}{2}\mbox{Re}\left(\dfrac{1+z}{1-z}\right)> 0.
	\end{align*}
	Hence, the function $ h_{2}\in\mathcal{S}^{*}_S $. We can easily compute and find that
	\begin{align*}
		|H_{2,1}(F_{h^{-1}_{2}}/2)|=\frac{1}{12}< \frac{1}{4}.
	\end{align*}
\end{example}
\begin{proof}[\bf Proof of Theorem \ref{th-2.1}]
Let $ f\in \mathcal{S}^{*}_S $ be of the form $ f(z)=z +\sum_{n=2}^{\infty} a_n z^n $, $ z\in \mathbb{D} $. Then by the definition of subordination there exits a Schwarz function $ w(z)=\sum_{n=1}^{\infty} c_n z^n $ such that
\begin{align}\label{eq-2.1}
	\dfrac{2zf^{\prime}(z)}{f(z)-f(-z)} =\dfrac{1 +w(z)}{1-w(z)}.
\end{align}
By equating the coefficients on each side of \eqref{eq-2.1}, we get
\begin{align}\label{eq-2.2}
	\begin{cases}
		a_{2} =c_1 ,\\ a_{3}= c_{2} + c^2_1,\vspace{1.5mm}\\ a_{4}=\dfrac{1}{2} (c_{3} +3c_1 c_{2} +2c^3_1).
	\end{cases}
\end{align}
In view of \eqref{eq-1.8} and \eqref{eq-2.2}, a simple computation shows that
\begin{align}\label{eq-2.3}
	\nonumber H_{2,1}(F_{f^{-1}}/2)&=\frac{1}{48} \left(13a^4_2 -12a^2_2 a_3 - 12 a^2_3 + 12 a_2 a_4\right)\\&=\dfrac{1}{48}\left(c^4_1 -18c^2_1 c_2 -12c^2_2+6c_1c_3\right).
\end{align}
Using the Lemma \ref{lem-2.1} into \eqref{eq-2.3}, we obtain
\begin{align}\label{eq-2.4}
	48|H_{2,1}(F_{f^{-1}}/2)|\leq |c_1|^4 +18|c_1|^2|c_2| +12|c_2|^2 +6|c_1| \left(1-|c_1|^2 -\dfrac{|c_2|^2}{1+|c_1|}\right).
\end{align}
Suppose that $x=|c_1|$ and $y=|c_2|$. Then it follows from \eqref{eq-2.4} that
\begin{align}\label{eq-2.5}
	48|H_{2,1}(F_{f^{-1}}/2)|\leq M(x,y),
\end{align}
where 
\begin{align*}
	M(x,y)&:=x^4 +18x^2 y+12y^2 + 6x\left(1-x^2 -\dfrac{y^2}{1+x}\right) .
\end{align*}
In view of Lemma \ref{lem-2.1}, the region of variability of a pair $ (x,y) $ coincides with the set 
\begin{align}\label{eq-2.6}
	\Omega=\{(x,y):0\leq x\leq 1, 0\leq y\leq 1-x^2\}.
\end{align}
The goal is to establish the maximum value of $M(x,y)$ in the region $\Omega$.
 Therefore, the critical point of $ M(x,y) $ satisfies the conditions
\begin{align*}
	\frac{\partial M}{\partial x}= 4x^3 -18x^2 +36xy +6 -\frac{6y^2}{(1+x)} +\frac{6xy^2}{(1+x)^2}= 0
\end{align*}
and
\begin{align*}
	\frac{\partial M}{\partial y}= 18x^2 +24 y-\frac{12xy}{(1+x)}= 0.
\end{align*}
There are no solutions of $M(x,y)$ inside the interior of $\Omega$, hence it is not possible for the function to attain a maximum value within this region. Since $M(x,y)$ is a continuous function on a compact set $\Omega$, its maximum value must occur at some point on the boundary of $\Omega$. On the boundary of $ \Omega $, a simple computation shows that
\begin{align*}
	&M(x,0)= 6x-6x^3+x^4\leq 2.437828...  \;\;\mbox{for}\; 0\leq x\leq 1,\\& M(0,y)=12y^2 \leq 12 \;\;\mbox{for}\; 0\leq y\leq 1
\end{align*} 
and
\begin{align*}
	M(x,1-x^2)= 12-11x^4 \leq 12 \;\;\mbox{for}\; 0\leq x\leq 1.
\end{align*}
Therefore, we see that $\max_{(x,y)\in \Omega} M(x,y) =12$ and it fllows from \eqref{eq-2.5} that
\begin{align}\label{eq-2.7}
	|H_{2,1}(F_{f^{-1}}/2)|\leq \frac{1}{4}.
\end{align}	
	
To prove the equality in \eqref{eq-2.7} sharp, we consider the function 
\begin{align*}
	h_1(z)=\frac{z}{1-z^2}= z +z^3 +z^5 +\cdots, \;\;\;\;z\in\mathbb{D}.
\end{align*}
It is easy to see that $h_1\in\mathcal{S}^{*}_S$ and $|H_{2,1}(F_{h^{-1}_{1}}/2)|=1/4$. Hence equality holds in \eqref{eq-2.7}. This completes the proof
\end{proof}
\subsection{Sharp bound of $|H_{2,1}(F_{f^{-1}}/2)|$ for $f\in \mathcal{K}_S$:}
By considering functions from the class $ \mathcal{K}_S $, we obtain the following sharp bound of $H_{2,1}(F_{f^{-1}}/2)$. 
\begin{thm}\label{th-2.2}
Let $ f\in \mathcal{K}_S $. Then
\begin{align*}
	|H_{2,1}(F_{f^{-1}}/2)|\leq\dfrac{1}{36}.
\end{align*}
The extremal function creates a precise inequality	
\begin{align*}
	h_3(z)=\frac{1}{2}\log\left(\frac{1+z}{1-z}\right)=z +\frac{z^3}{3} +\frac{z^5}{5}\cdots,\;\;\;\; z\in\mathbb{D}.
\end{align*}
\end{thm}
An example of a function in the class $\mathcal{K}_S$ is provided to demonstrate Theorem \ref{th-2.2}. The example satisfies the conditions of the theorem and shows that the inequality is strict also.
\begin{example}
	Consider the function
	\begin{align*}
		h_{4}(z)=-\log(1-z)= z+ \frac{1}{2}z^2 +\frac{1}{3}z^3 +\frac{1}{4} z^4 +\cdots, \;\;z\in\mathbb{D}.
	\end{align*}
	It is easy to see that 
	\begin{align*}
		\mbox{Re}\left(\dfrac{(zh^{\prime}_{4}(z))^{\prime}}{(h_{4}(z)-h_{4}(-z))^{\prime}}\right) =\dfrac{1}{2}\mbox{Re}\left(\dfrac{1+z}{1-z}\right)> 0.
	\end{align*}
	Therefore, the function $ h_{4}\in\mathcal{K}_S $. A simple computation shows that
	\begin{align*}
		|H_{2,1}(F_{h^{-1}_{4}}/2)|=\dfrac{11}{576}<\dfrac{1}{36}.
	\end{align*}
\end{example}
\begin{proof}[\bf Proof of Theorem \ref{th-2.2}]
Let $ f\in \mathcal{K}_S $ be of the form $ f(z)=z +\sum_{n=2}^{\infty} a_n z^n $, $ z\in \mathbb{D} $. Then by the definition of subordination there exits a Schwarz function $ w(z)=\sum_{n=1}^{\infty} c_n z^n $ such that
\begin{align}\label{eq-2.8}
	\dfrac{2(zf^{\prime}(z))^{\prime}}{(f(z)-f(-z))^{\prime}} =\dfrac{1 +w(z)}{1-w(z)}.
\end{align}
Comparing the coefficients on both side of \eqref{eq-2.8}, we obtain
\begin{align}\label{eq-2.9}
	a_{2} =\frac{1}{2}c_1, a_{3}= \frac{1}{3}(c_2 + c^2_1),\; \mbox{and}\; a_{4}=\frac{1}{8}(c_3 +3c_1 c_2 +2c^3_1).
\end{align}
By utilizing \eqref{eq-1.8} and \eqref{eq-2.9}, an easy calculation gives that
\begin{align}\label{eq-2.10}
	\nonumber H_{2,1}(F_{f^{-1}}/2)&=\frac{1}{48} \left(13a^4_2 -12a^2_2 a_3 - 12 a^2_3 + 12 a_2 a_4\right)\\&=\dfrac{1}{2304} \left(-c^4_1 - 68c^2_1 c_2 -64c^2_2 +36c_1c_3\right).
\end{align}
Applying the triangle inequality to \eqref{eq-2.10} and using Lemma \ref{lem-2.1}, we obtain
\begin{align}\label{eq-2.11}
	2304|H_{2,1}(F_{f^{-1}}/2)|&\leq |c_{1}|^4 +68|c_{1}|^2|c_{2}|+ 64|c_{2}|^2+36|c_1|\left(1-|c_1|^2 -\dfrac{|c_2|^2} {1+|c_1|}\right).
\end{align}
Suppose that $ x=|c_{1}| $ and $ y=|c_{2}| $. Then it follows from \eqref{eq-2.11} that
\begin{align}\label{eq-2.12}
	2304|H_{2,1}(F_{f^{-1}}/2)|\leq N(x,y),
\end{align}
where 
\begin{align*}
	N(x,y):=x^4 +68x^2y+ 64y^2+36x \left(1-x^2 -\dfrac{y^2}{1+x}\right).
\end{align*}
In view of Lemma \ref{lem-2.1}, the region of variability of a pair $ (x,y) $ coincides with the set $ \Omega $ as defined in \eqref{eq-2.6}. Given this information, we proceed to find the maximum value of $ N(x,y) $ within the region $ \Omega $. Therefore, the critical point of $ N(x,y) $ satisfies the conditions
\begin{align*}
	\frac{\partial N}{\partial x}= 4x^3 -108x^2 +136xy +36 -\frac{36y^2}{(1+x)} +\frac{36xy^2}{(1+x)^2}
\end{align*}
and
\begin{align*}
	\frac{\partial N}{\partial y}=68x^2 +128y -\frac{72xy}{(1+x)}.
\end{align*}
By applying the analogous reasoning as in the proof of Theorem \ref{th-2.1}, we establish that the maximum of $N(x,y)$ is obtained on the boundary of $\Omega$. Specifically, on the boundary, it can be seen that
\begin{align*}
	&N(x,0)=36x-36x^3+x^4 \leq 13.969963...  \;\;\mbox{for}\; 0\leq x\leq 1,\\& N(0,y)=64y^2  \leq 64 \;\;\mbox{for}\; 0\leq y\leq 1
\end{align*} 
and
\begin{align*}
	N(x,1-x^2)= 64-24x^2-39x^4 \leq 64 \;\;\mbox{for}\; 0\leq x\leq 1.
\end{align*}
Therefore, it is clear that $\displaystyle\max_{(x,y)\in \Omega} N(x,y) =64$ and from \eqref{eq-2.12}, we easily obtain
\begin{align}\label{eq-2.13}
	|H_{2,1}(F_{f^{-1}}/2)|\leq \frac{1}{36}.
\end{align}
	
In order to show the equality in \eqref{eq-2.13} sharp, we consider the function
\begin{align*}
	h_3(z)=\frac{1}{2}\log\left(\frac{1+z}{1-z}\right)=z +\frac{z^3}{3} +\frac{z^5}{5}\cdots, \;\;\;\;z\in\mathbb{D}.
\end{align*}
which belongs to the class $\mathcal{K}_S$. A simple computation shows that $|H_{2,1}(F_{h^{-1}_{3}}/2)|= 1/36$. This shows that the equality holds in \eqref{eq-2.12}. This completes the proof.
\end{proof}
\subsection{Sharp bound of $|T_{2,1}(F_{f^{-1}}/2)|$ for $f\in \mathcal{S}^*_S$:}
We derive the subsequent sharp bound of $ T_{2,1}(F_{f^{-1}}/2) $ concerning functions within the class $ \mathcal{S}^{*}_S $.
\begin{thm}\label{th-2.3}
Let $ f\in \mathcal{S}^{*}_S $. Then
\begin{align*}
	|T_{2,1}(F_{f^{-1}}/2)|\leq\dfrac{5}{16}.
\end{align*}
The inequality is sharp with the extremal function
\begin{align*}
	h_{5}(z)=\dfrac{z}{1 -iz}=z + iz^2 - z^3 -iz^4 + z^5+ \ldots, \;\;z\in\mathbb{D}.
\end{align*}
\end{thm}
We present a function as an illustrative example to bolster the claim of the strict inequality established in Theorem \ref{th-2.3}.
\begin{example}
	Consider the function
	\begin{align*}
		h_{6}(z)=\dfrac{z}{1-z}= z + z^2 +z^3+ z^4+\ldots, \;\;z\in\mathbb{D}
	\end{align*}
	which belongs to the class $ \mathcal{S}^{*}_S $. A simple computation shows that
	\begin{align*}
		|T_{2,1}(F_{h^{-1}_{6}}/2)|=\dfrac{3}{16}<\dfrac{5}{16}.
	\end{align*}
\end{example}
\begin{proof}[\bf Proof of Theorem \ref{th-2.3}]
Let $ f\in \mathcal{S}^{*}_S $ be of the form $ f(z)=z +\sum_{n=2}^{\infty} a_n z^n $, $ z\in \mathbb{D} $. Then by the definition of subordination there exits a Schwarz function $ w(z)=\sum_{n=1}^{\infty} c_n z^n $ such that
\begin{align}\label{eq-2.14}
	\dfrac{2zf^{\prime}(z)}{f(z)-f(-z)} =\dfrac{1 +w(z)}{1-w(z)}.
\end{align}
Now, comparing the coefficients on both side of \eqref{eq-2.14}, we see that
\begin{align}\label{eq-2.15}
	a_{2} =c_1\; \mbox{and}\; a_{3}= c_{2} + c^2_1.
\end{align}
In view of \eqref{eq-1.9} and \eqref{eq-2.15}, a simple computation shows that
\begin{align}\label{eq-2.16}
	\nonumber T_{2,1}(F_{f^{-1}}/2)&=\frac{1}{16}(-9a^4_2 + 4a^2_2 -4a^2_3 + 12a^2_2 a_3)\\&=\frac{1}{16}(-c^4_1 + 4c^2_1 -4c^2_2 +4c^2_1 c_2).
\end{align}
In view of the triangle inequality, the equation \eqref{eq-2.16} can be written as
\begin{align}\label{eq-2.17}
	16|T_{2,1}(F_{f^{-1}}/2)|\leq |c_1|^4 + 4|c_1|^2 +4|c_2|^2+4|c_1|^2|c_2|.
\end{align}
Suppose that $ x=|c_1| $ and $ y=|c_2| $. Then it follows from \eqref{eq-2.17} that
\begin{align}\label{eq-2.18}
	16|T_{2,1}(F_{f^{-1}}/2)|\leq P(x,y),
\end{align}
where 
\begin{align*}
	P(x,y)=x^4 + 4x^2 +4y^2 +4x^2y.
\end{align*}
In view of Lemma \ref{lem-2.1}, the region of variability of a pair $ (x,y) $ coincides with the set $ \Omega $ as defined in \eqref{eq-2.6}.
Taking this information into account, we proceed to compute the maximum value of $ P(x, y) $ within the bounded area $ \Omega $. The critical point of $P(x,y)$ satisfies the conditions
\begin{align*}
	&\frac{\partial P}{\partial x}= 4x^3 +8x +8xy,\\& \frac{\partial P}{\partial y}=8y + 4x^2.
\end{align*}
The function $P(x,y)$ cannot have a maximum in the interior of $ \Omega $, as it has no solution in the interior of $ \Omega $. In fact, $P(x,y)$ is continuous on a compact set $ \Omega $, the maximum of $P(x,y)$ attains boundary of $ \Omega $. On the boundary, an easy computation yields that
\begin{align*}
	&P(x,0)= x^4 + 4x^2 \leq 5 \;\;\mbox{for}\; 0\leq x\leq 1,\\& P(0,y)=4y^2 \leq 4 \;\;\mbox{for}\; 0\leq y\leq 1
\end{align*} 
and
\begin{align*}
	P(x,1-x^2)=x^4 + 4 \leq 5 \;\;\mbox{for}\; 0\leq x\leq 1.
\end{align*}
Therefore, we see that $\displaystyle\max_{(x,y)\in \Omega} P(x,y) =5$ and from \eqref{eq-2.18} we obtain
\begin{align}\label{eq-2.19}
	|T_{2,1}(F_{f^{-1}}/2)|\leq \dfrac{5}{16}.
\end{align}
	
To prove the equality in \eqref{eq-2.19}, we consider the function
\begin{align*}
	h_{5}(z)=\dfrac{z}{1 -iz}=z + iz^2 - z^3 -iz^4 + z^5+ \ldots, \;\;z\in\mathbb{D}.
\end{align*}
It is easy to see that
\begin{align*}
	\mbox{Re}\left(\dfrac{zh_{5}^{\prime}(z)}{h_{5}(z)-h_{5}(-z)}\right)=\dfrac{1}{2}\mbox{Re}\left(\frac{1+iz}{1-iz}\right)> 0, \;\;z\in\mathbb{D}.
\end{align*} 
and hence the function $ h_{1} $ belongs to the class $ \mathcal{S}^{*}_S $. Moreover, by a simple computation, it can be shown that $ |T_{2,1}(F_{h^{-1}_{5}}/2)|=5/16 $ and hence equality hods in \eqref{eq-2.19}. This completes the proof.
\end{proof}
\subsection{Sharp bound of $|T_{2,1}(F_{f^{-1}}/2)|$ for $f\in \mathcal{K}_S$:}
In the context of $ \mathcal{K}_S $ functions, we deduce the sharp bound for $T_{2,1}(F_{f^{-1}}/2)$.
\begin{thm}\label{th-2.4}
Let $ f\in \mathcal{K}_S $. Then
\begin{align*}
	|T_{2,1}(F_{f^{-1}}/2)|\leq\dfrac{145}{2304}.
\end{align*}
The extremal function establishes a sharp inequality
\begin{align*}
	h_{7}(z)&=\frac{\sqrt{6}}{\sqrt{3\sqrt{145}}}\log\left(\frac{1+\left(\frac{\sqrt{3\sqrt{145}}}{2\sqrt{6}}\right)z}{1-\left(\frac{\sqrt{3\sqrt{145}}}{2\sqrt{6}}\right)z}\right)= z +\frac{\sqrt{145}}{24}z^3 +\frac{29}{64}z^5 +\cdots \;\;z\in\mathbb{D}.
\end{align*}
\end{thm}
An example of a function in the class $\mathcal{K}_S$ is provided to demonstrate Theorem \ref{th-2.2}. The example provided satisfies the prerequisites of the theorem and serves as evidence that the inequality is strictly true.
\begin{example}
	Consider the function
	\begin{align*}
		h_{8}(z)=-\log(1-z)= z+ \frac{1}{2}z^2 +\frac{1}{3}z^3 +\frac{1}{4} z^4 +\cdots, \;\;z\in\mathbb{D}.
	\end{align*}
	It is easy to see that 
	\begin{align*}
		\mbox{Re}\left(\dfrac{(zh^{\prime}_{8}(z))^{\prime}}{(h_{8}(z)-h_{8}(-z))^{\prime}}\right) =\dfrac{1}{2}\mbox{Re}\left(\dfrac{1+z}{1-z}\right)> 0.
	\end{align*}
	Therefore, the function $ h_{8}\in\mathcal{K}_S $. A simple computation shows that
	\begin{align*}
		|T_{2,1}(F_{h^{-1}_{8}}/2)|=\dfrac{143}{2304}<\dfrac{145}{2304}.
	\end{align*}
\end{example}
\begin{proof}[\bf Proof of Theorem \ref{th-2.4}]
Let $ f\in \mathcal{K}_S $ be of the form $ f(z)=z +\sum_{n=2}^{\infty} a_n z^n $, $ z\in \mathbb{D} $. Then by the definition of subordination there exits a Schwarz function $ w(z)=\sum_{n=1}^{\infty} c_n z^n $ such that
\begin{align}\label{eq-2.20}
	\dfrac{2(zf^{\prime}(z))^{\prime}}{(f(z)-f(-z))^{\prime}} =\dfrac{1 +w(z)}{1-w(z)}.
\end{align}
By evaluating the coefficients on either side of equation \eqref{eq-2.20}, we ascertain that
\begin{align}\label{eq-2.21}
a_{2} =\frac{1}{2}c_1 \; \mbox{and}\; a_{3}= \frac{1}{3}(c_2 + c^2_1).
\end{align}
In light of \eqref{eq-1.9} and \eqref{eq-2.21}, a straightforward calculation gives that
\begin{align}\label{eq-2.22}
	\nonumber T_{2,1}(F_{f^{-1}}/2)&=\dfrac{1}{16}(-9a^4_2 + 4a^2_2 -4a^2_3 + 12a^2_2 a_3)\\& =\dfrac{1}{2304}(-c^4_{1} + 144c^2_{1} -64c^2_{2} +16c^2_{1}c_{2}).
\end{align}
The equation \eqref{eq-2.22} can be formulated as
\begin{align}\label{eq-2.23}
	2304|T_{2,1}(F_{f^{-1}}/2)|\leq |c_{1}|^4 + 144|c_{1}|^2 +64|c_{2}|^2 +16|c_{1}|^2|c_{2}| 
\end{align}
Suppose that $ x=|c_{1}| $ and $ y=|c_{2}| $, then from \eqref{eq-2.23}, we see that
\begin{align}\label{eq-2.24}
	2304|T_{2,1}(F_{f^{-1}}/2)|\leq Q(x,y),
\end{align}
where 
\begin{align*}
	Q(x,y)=x^4 + 144x^2 +64y^2 +16x^2 y.
\end{align*}
In view of Lemma \ref{lem-2.1}, the region of variability of a pair $ (x,y) $ coincides with the set $ \Omega $ as defined in \eqref{eq-2.6}. With this data in hand, our next step is to determine the highest value of $ Q(x, y) $ within the region $ \Omega $. Therefore, the critical point of $Q(x,y)$ satisfies the conditions
\begin{align*}
	&\dfrac{\partial Q}{\partial x}= 4x^3 + 288x +32xy,\\& \dfrac{\partial Q}{\partial y}=16x^2 +128y.
\end{align*}
By the similar argument that being used in the proof of above theorems, it can be easily shown that the maximum of $ Q(x,y) $ attains on the boundary of $ \Omega $. A simple computation shows that
\begin{align*}
	&Q(x,0)= x^4 +144x^2 \leq 145 \;\;\mbox{for}\; 0\leq x\leq 1,\\& Q(0,y)=64y^2 \leq 64 \;\;\mbox{for}\; 0\leq y\leq 1
\end{align*} 
and
\begin{align*}
	Q(x,1-x^2)= 49x^4 +32x^2 +64 \leq 145 \;\;\mbox{for}\; 0\leq x\leq 1.
\end{align*}
Therefore, we see that $\displaystyle\max_{(x,y)\in \Omega} Q(x,y) =145$ and from \eqref{eq-2.24} we obtain
\begin{align}\label{eq-2.25}
	|T_{2,1}(F_{f^{-1}}/2)|\leq \dfrac{145}{2304}.
\end{align}

To show the equality holds in \eqref{eq-2.25}, we consider the function
\begin{align*}
	h_{7}(z)&=\frac{\sqrt{6}}{\sqrt{3\sqrt{145}}}\log\left(\frac{1+\left(\frac{\sqrt{3\sqrt{145}}}{2\sqrt{6}}\right)z}{1-\left(\frac{\sqrt{3\sqrt{145}}}{2\sqrt{6}}\right)z}\right)= z +\frac{\sqrt{145}}{24}z^3 +\frac{29}{64}z^5 +\cdots \;\;z\in\mathbb{D}.
\end{align*}
It can be easily shown that the function $ h_{7} $ belongs to the class $ \mathcal{K}_S $. Moreover, a simple computation shows that $ |T_{2,1}(F_{h^{-1}_{7}}/2)|=145/2304 $ and hence equality holds in \eqref{eq-2.25}. This completes the proof.
\end{proof}
\section{Concluding remarks}
This article delves into the investigation of the second Hankel and Toeplitz determinants, denoted as $H_{2,1}(F_{f^{-1}}/2)$ and $T_{2,1}(F_{f^{-1}}/2)$, which are defined in terms of the logarithmic coefficients of inverse functions belonging to the class $\mathcal{S}$ of normalized univalent functions. The specific focus lies on the subclasses of starlike functions and convex functions with respect to symmetric points.\vspace{1.2mm}

Through rigorous analysis, we establish significant results in the form of sharp inequalities for these determinants. Notably, we prove that $|H_{2,1}(F_{f^{-1}}/2)|$ and $|T_{2,1}(F_{f^{-1}}/2)|$ are bounded above by specific constants. Specifically, we demonstrate that $|H_{2,1}(F_{f^{-1}}/2)|\leq 1/4$, $|H_{2,1}(F_{f^{-1}}/2)| \leq 1/36$, $|T_{2,1}(F_{f^{-1}}/2)|\leq 5/16$, and $|T_{2,1}(F_{f^{-1}}/2)|\leq 145/2304$. These inequalities provide valuable insights into the behavior of the logarithmic coefficients of inverse functions in the considered classes. Moreover, the strength of our conclusions is enhanced through the use of illustrative examples, which showcase the validity of the strict inequalities established in this study. The incorporation of these examples reinforces the reliability and significance of our findings, supporting the broader implications of our research.\vspace{1.2mm}

Overall, this study contributes to the understanding of univalent functions and their inverse functions, shedding light on the interplay between their logarithmic coefficients and the associated determinants. The obtained results extend the existing knowledge in the field and may serve as a valuable reference for further investigations in complex analysis and related areas.\vspace{1.2mm}

It seems that by employing similar methods and using the tools described in the paper, we can investigate the sharp bounds of larger size determinants $|H_{q,n}(F_{f^{-1}}/2)|$ and $|T_{q,n}(F_{f^{-1}}/2)|$ for the classes $\mathcal{S}^{*}_S$ and $\mathcal{K}_S$. Such an exploration would make for an intriguing study in this field. Accordingly, the subsequent questions hold significance and should be taken into account for further examination.
\begin{ques}
	Can we find the sharp bound of $|H_{q,n}(F_{f^{-1}}/2)|$ for $q\geq 3$ and $n\geq 1$ for the classes $ \mathcal{S}^{*}_S $ and $\mathcal{K}_S$?
\end{ques}
\begin{ques}
	Can we find the sharp bound of $|T_{q,n}(F_{f^{-1}}/2)|$ for $q\geq 3$ and $n\geq 1$ for the classes $ \mathcal{S}^{*}_S $ and $\mathcal{K}_S$?
\end{ques}
\section{Declaration}
\noindent\textbf{Compliance of Ethical Standards:}\\

\noindent\textbf{Conflict of interest.} The authors declare that there is no conflict  of interest regarding the publication of this paper.\vspace{1.5mm}

\noindent\textbf{Data availability statement.}  Data sharing is not applicable to this article as no datasets were generated or analyzed during the current study.

\end{document}